\documentclass{article}
\catcode`\@=11
\def\demo{\noindent{\bf Proof .-}}
\def\section{\@startsection {section}{1}{\z@}{-3.5ex plus -1ex
minus-.2ex}{2.3ex plus .2ex}{\normalsize\bf}}
\def\subsection{\@startsection {subsection}{1}{\z@}{-3.5ex plus -1ex
minus-.2ex}{2.3ex plus .2ex}{\normalsize\bf}}
\newtheorem{theorem}{Theorem}
\newtheorem{lemma}{Lemma}
\newtheorem{proposition}{Proposition}
\newtheorem{definition}{Definition}
\newtheorem{example}{Example}
\newtheorem{remark}{Remark}

\newtheorem{corollary}{Corollary}
\def\lcm{{\rm lcm}\,}

\begin{document}
\begin{center}
{\LARGE\bf{On the arithmetical rank of the edge ideals of forests}}
\end{center}
\vskip.5truecm
\begin{center}
{Margherita Barile\footnote{Partially supported by the Italian Ministry of University and Research.}\\ 
Dipartimento di Matematica, Universit\`{a} degli Studi di Bari\\
Via E. Orabona 4, 70125 Bari\\ Italy\\
e-mail: barile@dm.uniba.it}
\end{center}
\vskip1truecm
\noindent
{\small
{\bf Abstract} We show that for the edge ideals of a certain class of forests, the arithmetical rank equals the projective dimension.}

\section*{Introduction} Given a polynomial ring $R$ over a field, and a graph $G$ having the set of indeterminates as its vertex set, one can associate with $G$ a monomial ideal of $R$: this ideal is generated by the products of the vertices of each edge of $G$, and is hence generated by squarefree quadratic monomials. It is called the {\it edge ideal} of $G$, and was first introduced by Simis, Vasconcelos and Villarreal \cite{SVV}. It is interesting to derive algebraic invariants of this ideal from the combinatorial properties of the graph $G$. An extensive treatment of this kind is contained in the PhD thesis of Jacques \cite{J}, where the modules of the free resolutions of edge ideals are completely determined by means of a recursive construction.  In a recent work \cite{CN}, Corso and Nagel give a closed formula for the Betti numbers of the edge ideals of a special class of bipartite graphs,  the so-called {\it Ferrers graphs}. The same problem has been studied for the more general case of hypergraphs by Ha and van Tuyl in \cite{THvT}. 
  As an application of the results in \cite{CN}, in \cite{B2} it is shown that for the edge ideal of every Ferrers graph,  the arithmetical rank (ara), i.e., the least number of elements of $R$ which generate the ideal up to radical, is equal to the projective dimension (pd), i.e., to the length of every minimal free resolution of the quotient of $R$ with respect to the ideal.\newline In the present paper the same property is studied for the edge ideals of acyclic graphs, the so-called {\it forests}; from \cite{J} we know that in this case the projective dimension does not depend on the ground field (see also \cite{JK}).  We conjecture that for every edge ideal of a forest, the arithmetical rank equals the projective dimension: we, however, cannot prove this result in general. Using the recursive construction from \cite{J}, we can show the claim for a special class of forests, which we call {\it stretched}.  For these forests, the minimum number of elements generating the edge ideal up to radical can be given in form of a so-called {\it tree-like system}: this is a sequence of monomials and sums of two monomials which reflects the combinatorial structure of the forest.  We also give some classes of forests for which the arithmetical rank coincides with the general upper bound which was determined in \cite{B} for squarefree monomial ideals. Furthermore, we show that this upper bound can be re-formulated in terms of the number of edges and the degrees of the vertices of the graph. In the last section we explicitly compute the minimal free resolutions of the edge ideals of  a special class of (non stretched) forests and show that they are all 2-linear. To this end, we apply the approach presented by Lyubeznik in \cite{L1}. 
\newline
 It is worth pointing out that the computation of the arithmetical rank of an ideal in a polynomial ring is, in general, an open problem. For squarefree monomial ideals, a lower bound is provided by the projective dimension, but  there are not many known examples where this is the actual value of the arithmetical rank; some of them were classified by Terai and can be found in \cite{Te} and  \cite{Te2}.

\section{On the arithmetical rank of monomial ideals}
Consider the polynomial ring $R=K[x_1,\dots, x_n]$, where $K$ is a field. We recall some results on the arithmetical rank of the ideals in $R$ that are generated by monomials. Since the arithmetical rank does not change when an ideal is replaced by its radical, we can restrict our study to ideals generated by squarefree monomials. \newline
A finite set of elements of $R$  which generate a given ideal up to radical can be constructed according to the following well-known criterion, which is due to Schmitt and Vogel.
\begin{lemma}\label{lemma}{\rm [\cite{SV}, p.249]} Let $P$ be a finite subset of elements of $R$. Let $P_0,\dots, P_r$ be subsets of $P$ such that
\begin{list}{}{}
\item[(i)] $\bigcup_{i=0}^rP_i=P$;
\item[(ii)] $P_0$ has exactly one element;
\item[(iii)] if $p$ and $p'$ are different elements of $P_i$ $(0<i\leq r)$ there is an integer $i'$ with $0\leq i'<i$ and an element in $P_{i'}$ which divides $pp'$.
\end{list}
\noindent
We set $q_i=\sum_{p\in P_i}p^{e(p)}$, where $e(p)\geq1$ are arbitrary integers. We will write $(P)$ for the ideal of $R$ generated by the elements of $P$.  Then we get
$$\sqrt{(P)}=\sqrt{(q_0,\dots,q_r)}.$$
\end{lemma}
\noindent
In \cite{B} the above result is used to provide a general upper bound for the arithmetical rank of any squarefree monomial ideal $I$ of $R$. 
 Let $M=\{f_1,\dots, f_s\}$  be the set of its minimal monomial generators. Set $\mu(I)=\vert M\vert=s$. 
  Let $I_1, \dots, I_r$ be the minimal primes of $I$, so that $I=\cap_{j=1}^rI_j$.
Moreover, for all $i=1,\dots, s$, define
$$M_i=\{f_j: x_i\mbox{ divides }f_j\}.$$
\noindent
 For all $j=1,\dots, r$, set
\begin{equation}\label{nuj}\nu_j=\max\{\vert M_i\vert: x_i\in I_j\},\end{equation}
and
\begin{equation}\label{nu}\nu=\nu(I)=\min\{\nu_j: j=1,\dots, r\}.\end{equation}
\noindent Finally, for all $j=1,\dots, r$, set
\begin{equation}\label{roj}\rho_j=\min\{\vert M_i\vert:x_i \mbox{ divides } f_j\},\end{equation}
\noindent
and 
\begin{equation}\label{ro}\rho=\rho(I)=\max\{\rho_j:j=1,\dots, s\}.\end{equation}
\noindent
In \cite{B}, Proposition 1,  it is shown that
\begin{equation}\label{main}{\rm ara}\,I\leq \mu(I)-\nu(I)+1.\end{equation}
\noindent
The above upper bound can be rewritten in a different way, which, as we will see later, emphasizes its combinatorial character. In fact we have the following identity.
\begin{proposition}\label{Proposition} $\nu(I)=\rho(I)$.
\end{proposition}
\noindent
\demo We prove the two inequalities. First we show that $\nu\geq\rho$. Without loss of generality we may assume that $\rho=\rho_1$. For all $j=1,\dots, r$, there is a variable $x_{k_j}\in I_j$ such that $x_{k_j}$ divides $f_1$. Now, by (\ref{nuj}) and (\ref{roj}), for all $j=1\,\dots, r$, we have that 
$$\nu_j\geq\vert M_{k_j}\vert\geq \rho_1=\rho.$$
\noindent
By (\ref{nu}) this implies that $\nu\geq\rho$, as required. Next we show that  $\nu\leq\rho$. For all $j=1,\dots, s$, let $x_{h_j}$ be a variable dividing $f_j$ and such that $\rho_j=\vert M_{h_j}\vert$. Then $I\subset (x_{h_1},\dots, x_{h_s})$. Hence we may assume that $I_1\subset (x_{h_1},\dots, x_{h_s})$. Then by (\ref{nu}), (\ref{nuj}) and (\ref{ro}),
$$\nu\leq\nu_1\leq\max\{\rho_j:j=1,\dots, s\}=\rho.$$
\noindent
This completes the proof of the claim.
\par\bigskip\noindent
From \cite{L}  we know that, for any squarefree monomial ideal $I$,  the following inequality holds:
\begin{equation}\label{inequalities1}
{\rm pd}\, I\leq\,{\rm ara}\,I.
\end{equation}
\noindent
Moreover,  in view of (\ref{main}) and Proposition \ref{Proposition}, we have
\begin{equation}\label{inequalities2}
{\rm ara}\,I\leq \mu(I)-\rho(I)+1.
\end{equation}
\noindent
\begin{remark}\label{remark}{\rm Suppose that all $f_j$ have degree 2. In this case $I$ has a natural combinatorial interpretation. It can be associated with the following graph $G$ on the vertex set $\{x_1,\dots, x_n\}$:  
$$G=\{\{x_i, x_j\}:x_ix_j\in I\}.$$
\noindent
Then $I$ is the so-called {\it edge ideal} of the graph $G$, and we will denote it  by $I(G)$. Given an edge $\alpha=\{x_i, x_j\}$ of $G$, the monomial $a=x_ix_j$ will be called the {\it edge monomial} of $\alpha$. Thus we have that $I(G)$ is generated by the edge monomials of all edges of $G$. Hence $\mu(I)=\vert G\vert$, and 
$$\rho(I)=\max_{\alpha\in G}\{\min\{\deg(v):v\in \alpha\}\}.$$
}\end{remark}
Our goal is to determine the arithmetical ranks of the edge ideals of a big class of forests, which we will introduce in Section 3, and for which, as we will show, equality always holds in (\ref{inequalities1}). In Section 4 we will also present some examples where equality holds in  (\ref{inequalities2}), too. In the next section we develop the crucial tool that will be needed for the proof of our main result. 
\section{On tree-like systems}
This section is entirely devoted to the properties of the notion that is introduced by the following definition; we will also present its relevance in the generation of the edge ideals of forests up to radical. 
\begin{definition}\label{definition1}{\rm
Let $a_0, a_1, \dots, a_r\in R$ be non zero pairwise distinct elements and $b_1, b_2,\dots, b_n\in R$ be such that the non zero $b_i$'s are pairwise distinct and also distinct from the $a_i$'s. Call $P$ the set of all $a_i$'s and all non zero $b_i$'s. Moreover, set $q_0=a_0$ and $q_i=a_i+b_i$ for all $i=1,\dots, r$. Suppose that for all $i=1,\dots, r$ there is an index $j<i$ such that $a_j$ divides $a_ib_i$. Then the sequence
$$\Sigma: q_0, q_1, \dots, q_r$$
\noindent
is called a {\it tree-like system}.  The number $r+1$ is called the {\it length} of $\Sigma$ and is denoted by 
$\lambda(\Sigma)$.  The set $P$ is called the {\it support} of $\Sigma$. We will also say that $\Sigma$ {\it starts} at (the {\it starting point}) $a_0$ and {\it ends} at $q_r$. If $b_i=0$, we will say that $q_i=a_i$ is an {\it isolated summand}.\newline
A subsequence of $\Sigma$ which is a tree-like system will be called a {\it subtree} of $\Sigma$. 
}
\end{definition}
\begin{remark}\label{remark'}{\rm In the sequel we will often need to obtain a single tree-like system from several tree-like systems. This will happen according to two basic constructions. 
Given two tree-like systems $\Sigma_1: q^{(1)}_0, q^{(1)}_1,\dots, q^{(1)}_{r_1}$ with support $P^{(1)}$ and $\Sigma_2: q^{(2)}_0, q^{(2)}_1,\dots, q^{(2)}_{r_2}$ with support $P^{(2)}$, if $P^{(1)}$ and $P^{(2)}$ are disjoint, then
 $$\Sigma_1, \Sigma_2: q^{(1)}_0, q^{(1)}_1,\dots, q^{(1)}_{r_1}, q^{(2)}_0, q^{(2)}_1,\dots, q^{(2)}_{r_2}$$
\noindent
  is a tree-like system with support $P^{(1)}\cup P^{(2)}$. This construction by juxtaposition obviously extends to any finite number of tree-like systems with pairwise disjoint supports.\newline 
If $\Sigma_1$ and $\Sigma_2$ are subtrees of the same tree-like system $\Sigma$, then we denote by $\Sigma_1\cup\Sigma_2$ the subsequence formed by the elements of $\Sigma$ which belong to $\Sigma_1$ or $\Sigma_2$. This is, of course, a tree-like system. This construction obviously extends to any set of subtrees of $\Sigma$. 
}
\end{remark}
With respect to the notation of Definition \ref{definition1} we also have the next result.
\begin{proposition}\label{remark1}  
$$\sqrt{(q_0,q_1,\dots,q_r)}=\sqrt{(P)}.$$
\end{proposition}
\demo Set $P_0=\{a_0\}$ and $P_i=\{a_i, b_i\}$ for all $i=1,\dots, r$. Then the assumption of Lemma \ref{lemma} is fulfilled. The claim follows.
\par\bigskip\noindent
According to Proposition \ref{remark1}, two tree-like systems with the same support generate the same ideal up to radical. This justifies the following
\begin{definition}{\rm Two tree-like systems with the same support are called {\it equivalent}. 
}
\end{definition}
\begin{remark}\label{remark0}{\rm In general, given a tree-like system $\Sigma: q_0, q_1, \dots, q_r$, there can be other arrangements of the elements $q_i$ which are still tree-like systems. For instance, whenever $\Sigma'$ is a subtree of $\Sigma$, then a tree-like system equivalent to $\Sigma$ can be obtained by placing $\Sigma'$ at the beginning of $\Sigma$ and then listing the remaining elements of $\Sigma$ in their original order. In this case we will say  that $\Sigma'$ is {\it pushed to the top} of $\Sigma$. In particular, pushing any subsequence of isolated summands to the top produces an equivalent tree-like system.} 
\end{remark}
\begin{remark}\label{remark01}{\rm Let $S$ be a subset of the tree-like system $\Sigma$, and let $Q$ be the set of summands of $S$. Let $\Sigma'$ be a tree-like system with support $Q$. Then a tree-like system equivalent to $\Sigma$ can be obtained by first omitting the elements of $S$ and then placing $\Sigma'$ before the residual subsequence of $\Sigma$.  In this case we will simply say that $S$ is {\it replaced} by $\Sigma'$ in $\Sigma$. Note that this construction in particular applies when $S=\Sigma'$ is a subtree of $\Sigma$. 
}
\end{remark}
\begin{remark}\label{remark2}{\rm Let $\Sigma: q_0,q_1\,\dots, q_r$ be a tree-like system whose support is the set of edge monomials of the forest $T$. Then, for all $i=1,\dots, r$, if $q_i$ is not an isolated summand, we have that $q_i=uv+wz$, where $\{u,v\}$ and $\{w, z\}$ are distinct edges of $T$. By Definition \ref{definition1}, there is an index $j<i$ such that $q_j$ contains an edge monomial (other than $uv$ and $wz$) dividing $uvwz$, i.e., a monomial which is the product of a factor of $uv$ and a factor of $wz$. Up to renaming the indeterminates we may assume that this product is $vw$. Then  $\alpha=\{u,v\}$, $\gamma=\{v,w\}$, $\beta=\{w, z\}$ are three consecutive edges of $T$. We will say that the edge $\gamma$ {\it lies between} $\alpha$ and $\beta$.  Note that $vw$ is the only edge monomial of $T$ dividing $uvwz$: in fact none of $\{u,w\}$, $\{u,z\}$ and $\{v, z\}$ can be an edge of $T$, since otherwise $T$ would contain a cycle. This implies that, in the present case, the index $j$ in Definition \ref{definition1} is uniquely determined by $i$. We will call $q_j$ the {\it precedessor} of $q_i$ in $\Sigma$. We will also say that $q_i$ is a {\it follower} of $q_j$.\newline
Finally note that, since $\alpha$ and $\beta$ are disjoint edges, the vertices $u, v, w, z$ are pairwise distinct, i.e., no vertex can appear twice in any element of a tree-like system. 
 }
\end{remark}
\begin{definition}\label{definition2}{\rm In the assumption of Definition \ref{definition1}, if, for all $i=1,\dots, r$, $b_i\ne 0$ and $j=i-1$, then the tree-like system $\Sigma$ is called {\it strict}. 
}
\end{definition}
\begin{remark}\label{remark2'}{\rm As an immediate consequence of Definition \ref{definition2}, the starting point is the only isolated summand of a strict tree-like system. Moreover, if the support of a strict tree-like system is the set of edge monomials of a forest, from Remark \ref{remark2} we deduce that every element (except the last one) of this strict tree-like system contains an edge monomial whose corresponding edge lies between two other edges. In particular, if this strict tree-like system has more than one element, it cannot have the edge monomial of a terminal edge as its isolated summand.
}
\end{remark}
The strict tree-like systems are the fundamental constituents of the theory we are developing here. In fact every tree-like system is the union of strict subtrees, as we show next.
\begin{lemma}\label{strict}
For every element of a tree-like system there is a strict subtree ending at this element. If the support of the tree-like system is the set of edge monomials of a forest, then this strict subtree is unique.  
\end{lemma}
\demo  We refer to the notation introduced in Definition \ref{definition1}. Let $q$ be an element of $\Sigma$. If $q$ is an isolated summand, then it forms a strict subtree by itself. According to Remark \ref{remark2'}, this is the only subtree of $\Sigma$ ending at $q$.    So assume that $q=q_i=a_i+b_i$, where $a_i, b_i$ are distinct non zero elements of $R$. Set $i_0=i$ and, for $k\geq 0$, if $q_{i_{k}}$ is not an isolated summand, recursively define $i_{k+1}$ as the index such that $q_{i_{k+1}}$ is the precedessor of $q_{i_k}$ in $\Sigma$. Then the indices $i_k$ form a strictly descending sequence of nonnegative integers. This can only have finitely many terms; hence the process must stop, i.e, we have that, for some $k'$, $q_{i_{k'}}$ is an isolated summand. Then $\Sigma': q_{i_{k'}}, q_{i_{k'-1}}, \dots, q_{i_0}$ is, by construction, a strict subtree of $\Sigma$ ending at $q=q_{i_0}$. If the support of $\Sigma$ is the set of edge monomials of a forest, then the uniqueness of $\Sigma'$ follows from the uniqueness of precedessors established in Remark \ref{remark2}.  \par\bigskip\noindent
We will say that two edges $\alpha$ and $\beta$ of a graph $T$ are {\it connected} if $T$ contains a sequence ({\it path}) of consecutive edges starting at $\alpha$ and ending at $\beta$. This defines an equivalence relation in the set of edges of $T$; the equivalence classes are the so-called {\it connected components} of $T$. Note that if $T$ is a forest, the aforementioned path is unique.
\begin{lemma}\label{path}
Let $\Sigma$ be a tree-like system whose support is the set of edge monomials of a forest $T$. Let $\alpha$ and $\beta$ be edges of $T$ such that the corresponding edge monomials $a$ and $b$ belong to the support of the same strict subtree $\Sigma'$ of $\Sigma$. Then 
there is a path of $T$ connecting $\alpha$ and $\beta$ (i.e., $\alpha$ and $\beta$ belong to the same connected component of $T$).
\end{lemma}
\demo According to Remark \ref{remark2} the claim is true if $a$ and $b$ are summands of the same element of $\Sigma'$, since in this case there is an edge lying between $\alpha$ and $\beta$. For the rest, it suffices to prove the claim in the case where $a$ and $b$ appear in two consecutive elements of $\Sigma'$, say $a$ in the precedessor of the element containing $b$. Assume that these elements are $q'_i=a+c$ and $q'_{i+1}=b+d$, for some $c,d\in R$. Then, by Remark \ref{remark2'}, $d\ne0$. Let $\delta$ be the corresponding edge of $T$. If $c\ne0$, let $\gamma$ be the corresponding edge of $T$. According to Remark \ref{remark2}, we then have that either $\alpha$ or $\gamma$ lies between $\beta$ and $\delta$. In the former case $\alpha$ and $\beta$ are consecutive edges, so that the claim is true. Let us consider the latter case. Then $q'_i$ is not an isolated summand, so there is an element $q'_{i-1}$ that is the precedessor of $q'_i$ in $\Sigma'$; it contains an edge monomial $e$ such that the corresponding edge $\varepsilon$  lies between $\alpha$ and $\gamma$. Then the path connecting $\alpha$ and $\beta$ is either $\alpha, \varepsilon, \gamma, \beta$ or $\alpha, \varepsilon, \beta$, and the claim is true in this case, too. 
\begin{corollary}\label{corollary_connected} Let $\Sigma$ be a tree-like system whose support is the set of edge monomials of the forest $T$, and let $C$ be a connected component of $T$. Moreover, let $\Sigma(C)$ be the set of all elements of $\Sigma$ that contain the edge monomial of an edge of $C$. Then $\Sigma(C)$ is a subtree of $\Sigma$ whose support is the set of edge monomials of $I(C)$. 
\end{corollary}
\demo Let $\alpha$ be an edge of $C$. By Lemma \ref{strict} there is a strict subtree $\Sigma_{\alpha}$ of $\Sigma$ whose support contains the edge monomial of $\alpha$. Moreover, by Lemma \ref{path}, all elements in the support of $\Sigma_{\alpha}$ correspond to edges of $C$. This proves that
$$\Sigma(C)=\bigcup_{\alpha\in C} \Sigma_{\alpha}$$
\noindent
is the required subtree of $\Sigma$; the union on the right-hand side is the one described in Remark \ref{remark'}. 
\par\bigskip\noindent
The next result presents an important combinatorial construction on tree-like systems which will play an important role in the proof of our main theorem.
\begin{lemma}\label{inversion}
Let $r\geq2$ be an integer and let $a_0, a_1,\dots, a_r, b_1,\dots, b_r\in R$ be pairwise distinct squarefree quadratic monomials such that, for all $i=1,\dots r$, $a_{i-1}$ divides $a_ib_i$. Then there is a strict tree-like system with support $\{a_0,a_1,\dots, a_r, b_1,\dots, b_r\}$ and starting point $a_{r-1}$.
\end{lemma}  
\demo
We proceed by induction on $r$. Consider the sets $S=\{a_0\}$, $L=\{a_1, a_2,\dots, a_r\}$ and $R=\{b_1, b_2,\dots, b_r\}$, whose elements are the starting point, the left summands and the right summands, respectively, of the strict tree-like system
$$\Sigma: q_0=a_0,\ q_1=a_1+b_1,\ \dots,\ q_r=a_r+b_r.$$
\noindent
First assume that $r=2$. Let $a_0=xy$, where $x,y$ are indeterminates. Since $a_0$ divides $a_1b_1$, up to renaming the indeterminates we may assume that $x$ divides $a_1$ and $y$ divides $b_1$. Then $a_1=xz$ for some indeterminate $z$ other than $y$ and $x$. Since, in turn, $a_1$ divides $a_2b_2$, up to interchanging $a_2$ and $b_2$ (which does not affect the assumption, since $a_2$ and $b_2$ have no followers), we may assume that $x$ divides $a_2$, and $z$ divides $b_2$. It follows that $a_0$ divides $a_2b_1$. Moreover, $a_1$ divides $a_0b_2$. Consequently, $\Sigma': q_0=a_1, q_1=a_0+b_2, q_2=a_2+b_1$ is the required strict tree-like system.  Now assume that $r>2$ and suppose the claim true for all smaller $r$. The induction basis applies to $S'=\{a_{r-2}\}$, $L'=\{a_{r-1}, a_r\}$, $R'=\{b_{r-1}, b_r\}$, since $a_{r-2}$ divides $a_{r-1}b_{r-1}$ and $a_{r-1}$ divides $a_rb_r$: thus the first part of the proof shows that (up to interchanging $a_r$ and $b_r$), 
$$\bar\Sigma: \bar q_0=a_{r-1},\ \bar q_1=a_{r-2}+b_r,\ \bar q_2=a_r+b_{r-1}$$
\noindent
is a tree-like system, where $a_{r-2}$ divides $a_rb_{r-1}$. By virtue of this latter statement, induction applies to $S''=\{a_0\}$, $L''=\{a_1, a_2,\dots, a_{r-2},a_r\}$ and $R''=\{b_1, b_2,\dots, b_{r-2}, b_{r-1}\},$ so that we have a tree-like system with support $S''\cup L''\cup R''$ and starting point $a_{r-2}$,
$\Sigma'': q''_0=a_{r-2},\ q''_1,\dots, q''_{r-1}.$ But then 
$$\Sigma''':q'''_0=a_{r-1},\ q'''_1=a_{r-2}+b_r, q'''_2=q''_1,\dots, q'''_r=q''_{r-1}$$
\noindent
is the required strict tree-like system. This completes the proof.
\par\smallskip\noindent
\begin{remark}\label{remark3}{\rm The claim of Lemma \ref{inversion} can be rephrased as follows: given a strict tree-like system whose last two elements are $q_{r-1}=a_{r-1}+b_{r-1},\ q_r=a_r+b_r$, where $a_{r-1}$ divides $a_rb_r$, we can construct an equivalent strict tree-like system with starting point $a_{r-1}$. One of the summands in the last but one element of the initial tree-like system is pushed to the first position; therefore we will refer to this transformation as a {\it tree-inversion}. It has the following graph-theoretical interpretation. Suppose that a tree $T$ is constructed by the following recursive procedure, which is performed $r$ times, for a fixed integer $r\geq2$:\newline
Step 1: Draw an edge $\alpha_0$. Set $i=0$.\newline
Step 2: Draw two edges $\alpha_{i+1}, \beta_{i+1}$, so that $\alpha_i$ lies between $\alpha_{i+1}$ and $\beta_{i+1}$.\newline
Step 3: Replace $i$ with $i+1$. If $i<r$, go to Step 2, else end.\par\smallskip\noindent
The tree-inversion lemma states that the same tree $T$ can be constructed with a similar procedure, starting at edge $\alpha_{r-1}$.
}
\end{remark}
\noindent
\section{The edge ideals of forests}
In the sequel, $T$ will be a forest with at least one edge. We consider the edge ideal $I(T)$ of $T$ in the polynomial ring $R=K[V]$, where $V$ is the vertex set of $T$.   For the proof of our main theorem we need some preliminary results on forests, which are due to Jacques  \cite{J} and Jacques and Katzman \cite{JK}.
\begin{proposition}\label{components}{\rm(\cite{J}, Proposition 2.2.8)} Let the graph $G$ be the disjoint union of the subgraphs $G_1,\dots, G_s$. Then 
$${\rm pd}\,I(G)=\sum_{i=1}^s{\rm pd}\,I(G_i).$$
\end{proposition}
From this we deduce the following result.
\begin{corollary}\label{corollary_components} In the assumption of Proposition \ref{components}, suppose that ara\,$I(G_i)=\,{\rm pd}\,I(G_i)$ for all $i=1,\dots, s$. Then 
$${\rm ara}\,I(G)=\,{\rm pd}\,I(G).$$
\end{corollary}
\demo For all $i=1,\dots, s$, set $a_i=\,{\rm ara}\,I(G_i)$ and let $q^{(i)}_1, q^{(i)}_2,\dots, q^{(i)}_{a_i}\in R$ be such that 
$\sqrt{(q^{(i)}_1, q^{(i)}_2,\dots, q^{(i)}_{a_i})}=I(G_i)$. Then 
$$\sqrt{\sum_{i=1}^s(q^{(i)}_1,q^{(i)}_2,\dots,q^{(i)}_{a_i}})=\sqrt{\sum_{i=1}^s\sqrt{(q^{(i)}_1,q^{(i)}_2,\dots,q^{(i)}_{a_i})}}=\sqrt{\sum_{i=1}^sI(G_i)}=I(G),$$
\noindent
whence
$${\rm ara}\, I(G)\leq \sum_{i=1}^sa_i=\sum_{i=1}^s{\rm pd}\,I(G_i)= {\rm pd}\, I(G)$$
\noindent
by Proposition \ref{components}. Since, on the other hand, by (\ref{inequalities1}),  
$${\rm pd}\,I(G)\leq {\rm ara}\,I(G),$$
\noindent
equality holds. This completes the proof.
\begin{proposition}\label{vertex}{\rm (\cite{JK}, Proposition 4.1)} Let $T$ be a forest. If $T$  contains a vertex of degree at least 2, then there exists a vertex $v$ with at least two neighbours, such that all but one of its neighbours have degree 1. 
\end{proposition}
\noindent
If $T$ has a vertex of degree at least 2, let $v$ be a vertex fulfilling the assumption of Proposition \ref{vertex}. Otherwise let $v$ be any vertex of $T$; in this case $v$ has one only neighbour whose degree is equal to 1 as well. Let $v_1,\dots, v_n$ be the neighbours of $v$, where $v_1,\dots, v_{n-1}$ have degree 1 (or $v_1$ has degree 1 if $n=1$). 
Let $T'$ be the subgraph of $T$ induced on $V\setminus\{v_1\}$ and let $T''$ be the subgraph of $T$ induced on $V\setminus\{v, v_1, \dots, v_n\}$. Note that $T'$ and $T''$ are forests.
\begin{proposition}\label{pd}{\rm (\cite{J}, Theorem 9.4.17)} We have:
$${\rm pd}\,I(T)=\max\{{\rm pd}\,I(T'), {\rm pd}\, I(T'')+n\}.$$
\end{proposition}
\noindent
Next we introduce the class of forests which will be the central object of study in this section. 
\begin{definition}{\rm A forest will be called {\it stretched} if every edge of it has one vertex of degree at most 2.
}
\end{definition}
\noindent
Our main result is the following.
\begin{theorem}\label{theorem} Let $T$ be a stretched forest. Then ${\rm pd}\,I(T)=\,{\rm ara}\,I(T)$ and there is a tree-like system of length ara\,$I(T)$ whose support is the set of edge monomials of $I(T)$. 
\end{theorem}
\demo In the sequel, for the sake of simplicity, a tree-like system whose support is the set of edge monomials of a graph $G$ will be called a tree-like system for  $I(G)$.\newline
Note that it suffices to prove the claim for stretched forests without isolated vertices, since the isolated vertices do not appear in the generators of the edge ideal. 
 First assume that all vertices of $T$ have degree 1. Then $T$ consists of pairwise disjoint edges. The set of edge monomials thus forms a regular sequence of generators (which is also a tree-like system), so that $I(T)$ is a complete intersection, and pd\,$I(T)=\vert T\vert=\,{\rm ara}\,I(T)$. Hence the claim is true in this case. Next assume that $T$ has one vertex of degree greater than 1. By Proposition \ref{vertex} there is one vertex $v$ of $T$ with neighbours $v_1,\dots, v_n$, where $n\geq2$ and $v_1,\dots, v_{n-1}$ have degree 1. First assume that $v_n$ has degree 1, too. Then the subgraph $C$ of $T$ induced on $\{v, v_1,\dots, v_n\}$ is a connected component of $T$ whose edge monomials are $vv_1,\dots, vv_n$; it is the star-graph $S_n$ of Section 4 below. We thus have that ara\,$I(C)\leq n$. On the other hand, by \cite{J}, Theorem 5.4.11, pd\,$I(C)=n$, so that, by (\ref{inequalities1}), $n\leq\,{\rm ara}\,I(C)$. Thus  
ara\,$I(C)=\,{\rm pd}\,I(C)=n$,  and the edge monomials of $C$ give the required tree-like system for $I(C)$. Now, in view of Corollary \ref{corollary_components}, the claim is true if and only if it is true for every connected component of $T$. Thus it suffices to prove the claim in the case where $v_n$ has degree greater than 1, i.e., it has some neighbour other than $v$. Let $w_1,\dots, w_m$ $(m\geq 1)$ be the neighbours of $v_n$ other than $v$ and consider the subgraphs $T'$ and $T''$ defined above. For convenience of notation, in this proof $T'$ will denote the subgraph induced on $V\setminus\{v_{n-1}\}$. Induction applies to $T'$ and $T''$, since every subgraph of a stretched forest is stretched, too. Set $A'=\,{\rm pd}\,I(T')$ and $A''=\,{\rm pd}\,I(T'')$ and 
\begin{equation}\label{max}A=\max\{A', A''+n\}.\end{equation}
\noindent
Since $\{v, v_n\}$ is an edge of $T$, and $T$ is stretched, either $v$ has degree 2 (in which case $n=2$) or $v_n$ has degree 2 (in which case $m=1$). 
Consider the following auxiliary claim: there is a tree-like system $\Sigma: q_0, q_1, \dots, q_s$ of length at most $A$ for  $I(T)$  such that 
 
\begin{list}{}{}
\item{(i)}  if $n=2$, then 
\begin{eqnarray*}
q_0&=&vv_2,\\
q_1&=&vv_1+v_2w_1;
\end{eqnarray*}
\noindent
\item{(ii)} if $n\geq3$ and $m=1$, then, up to renaming the indices:
\begin{eqnarray*}
q_0&=&vv_n,\\ 
q_1&=&vv_1+v_nw_1,\\
q_2&=&vv_2,\\
&\vdots&\\
q_{n-1}&=&vv_{n-1}.
\end{eqnarray*}
\end{list}
\noindent
Since, by Proposition \ref{pd} and (\ref{inequalities1}), $A=\,{\rm pd}\, I(T)\leq\,{\rm ara}\, I(T)$, once that the auxiliary claim is proven, it will follow that ${\rm ara}\,I(T)=A$, so that $A$ will turn out to be the actual length of $\Sigma$.  Therefore, the auxiliary claim implies the theorem.\newline 
We show the auxiliary claim by induction on the number $N$ of vertices of $T$. First assume that $n=2$. 
\noindent
The minimum $N$ for $n=2$ and $m\geq1$ is $N=4$ and corresponds to the graph $T$ on the vertex set $V=\{v, v_1, v_2, w_1\}$ whose set of edge monomials is $\{vv_1, vv_2, v_2w_1\}$; this is the line graph $L_4$ which will be presented in Section 4. Here $m=1$. Note that in this case $T'=\{\{v, v_2\}, \{v_2, w_1\}\}$ and $T''=\emptyset$, so that $A'=2$, $A''=0$, and, consequently, $A=2$. Then $\Sigma: q_0=vv_2, q_1=vv_1+v_2w_1$ is a tree-like system for $I(T)$ that fulfills the auxiliary claim in case (i). Now assume that $n=2$, $N>4$ and that the claim of the theorem is fulfilled by all stretched forests with less than $N$ vertices. In particular we will assume that $A'=$\,ara\,$I(T')$ and $A''=$\,ara\,$I(T'')$, and that there are a tree-like system $\Sigma': q'_0,\dots, q'_{A'-1}$ for $I(T')$ and a tree-like system $\Sigma'': q''_0,\dots, q''_{A''-1}$ for $I(T'')$.  Recall that 
\begin{eqnarray}\label{T} I(T)&=&I(T')+(vv_1),\nonumber\\
I(T')&=&I(T'')+(vv_2,v_2w_1,\dots, v_2w_m).\end{eqnarray}
\noindent
We will deduce that the auxiliary claim is true for $T$, i.e., that there is a tree-like system $\Sigma: q_0, q_1,\dots, q_{A-1}$ for $I(T)$ for which  (i) is fulfilled. 
\noindent
We distinguish between several cases, depending on where the edge monomials $vv_2, v_2w_1, \dots, v_2w_m$ of $T'$ appear in $\Sigma'$.
\par\bigskip\noindent
\underline{Case 1}: The edge monomials $vv_2$ and $v_2w_i$, for some index $i$, $1\leq i\leq m$, both appear as isolated summands in $\Sigma'$. We may assume that $i=1$. According to Remark \ref{remark0}, up to rearranging the elements of $\Sigma'$, we may also assume that $q'_0=vv_2$ and $q'_1=v_2w_1$. Set $q_0=q'_0$, $q_1=vv_1+v_2w_1$, and $q_i=q'_i$ for $i=2,\dots, A'-1$. Then $\Sigma: q_0, q_1, \dots, q_{A'-1}$ is a tree-like system of length $A'\leq A$ for $I(T)$ fulfilling (i). Thus the auxiliary claim is true in Case 1.
\par\bigskip\noindent
\underline{Case 2}: The edge monomial $vv_2$ appears in $\Sigma'$ as an isolated summand, but none of $v_2w_i$ does. Up to rearrangement we may assume that $q'_0=vv_2$. For all $i=1,\dots, m$, let $q'_{k_i}$ be the element of $\Sigma'$ containing the edge monomial $v_2w_i$ and let $x_i$ and $y_i$ be indeterminates such that 
\begin{equation}\label{q} q'_{k_i}=v_2w_i+x_iy_i.\end{equation}
\noindent Consider the edges $\alpha_i=\{v_2, w_i\}$, $\beta_i=\{x_i, y_i\}$ of $T'$. By Remark \ref{remark2}, for all $i=1,\dots, m$, there is an edge $\gamma_i$ lying between $\alpha_i$ and $\beta_i$. For every index $i$, one of $x_i$ and $y_i$, say $x_i$, belongs to $\gamma_i$, and, similarly,  either $v_2$ or $w_i$ belongs to $\gamma_i$. If $v_2\in \gamma$, then $x_i$ is a neighbour of $v_2$ other than $w_i$; it is also distinct form $v$, since $v$ is a terminal vertex of $T'$ and, consequently, cannot belong to $\gamma_i$. Thus for some index $j_i\ne i$, we have that $x_i=w_{j_i}$, i.e., $\gamma_i=\{v_2, w_{j_i}\}$,  and $\beta_i=\{w_{j_i}, y_i\}$, whence
\begin{equation}\label{q1} q'_{k_i}=v_2w_i+w_{j_i}y_i.\end{equation}
\noindent
If $w_i\in\gamma$, then  $\gamma_i=\{w_i, x_i\}$ and $y_i$ is not a neighbour of $v_2$, because otherwise the vertices $v_2, w_i, x_i, y_i$ would form a cycle. Therefore 
\begin{equation}\label{q2} q'_{k_i}=v_2w_i+x_iy_i,\qquad\mbox{where $x_i\ne v_2$ is a neighbour of $w_i$ and $y_i\ne w_j$ for all $j=1,\dots, m$}.\end{equation}
\noindent
At this point we have to distinguish between two subcases.
\par\smallskip\noindent
\underline{Case 2.1}: For some index $i\in\{1,\dots, m\}$, $\gamma_i=\{v_2, w_{j_i}\}$, i.e., $q'_{k_i}$ is of the form (\ref{q1}).  We may assume that $i=1$, and $j_1=2$.   By Lemma \ref{strict} there is a unique strict subtree $\bar\Sigma$  of $\Sigma'$ that ends at $q'_{k_1}=v_2w_1+w_2y_1$. According to Remark \ref{remark2}, the precedessor of $q'_{k_1}$ in  $\bar\Sigma$ contains the edge monomial $c_1=v_2w_2$ of the edge $\gamma_1$, hence it is $q'_{k_2}=v_2w_2+x_2y_2$.  Thus Lemma \ref{inversion} can be applied to perform a tree-inversion on $\bar\Sigma$ with $a_r=v_2w_1$, $b_r=w_2y_1$, $a_{r-1}=v_2w_2$, $b_{r-1}=x_2y_2$ so as to produce a strict tree-like system $\tilde\Sigma$ equivalent to $\bar\Sigma$ and with starting point $\tilde q_1=v_2w_2$. By Remark \ref{remark2'}, since $\{v, v_2\}$ is a terminal edge of $T'$, $q'_0=vv_2$ cannot belong to $\bar\Sigma$, hence it does not appear in $\tilde\Sigma$ either. If we replace $\tilde\Sigma$ for $\bar\Sigma$ in $\Sigma'$, according to Remark \ref{remark01}, we thus obtain a tree-like system for $I(T')$ where $vv_2$ and $v_2w_2$ are isolated summands. This takes us back to Case 1. Hence the auxiliary claim is true in Case 2.1.  
\par\smallskip\noindent
\underline{Case 2.2}: For all indices $i=1,\dots, m$, $\gamma_i=\{w_i, x_i\}$, i.e., $q'_{k_i}$ is of the form (\ref{q2}).  Note that  the vertices $w_1,\dots, w_m$ belong to pairwise distinct connected components of $T''$: if there were a path connecting $w_i$ and $w_j$ in $T''$ for some distinct indices $i$ and $j$, then this path would not contain the vertex $v_2$, so that the edge $\{v_2, w_i\}$, this path and the edge $\{v_2, w_j\}$ would form a cycle. For all $i=1,\dots, m$, let $C_i$ be the connected component of $T''$ containing $w_i$, and let $C_{m+1}, \dots, C_s$ be the remaining connected components of $T''$ (if any exist). Note that $C_{m+1}, \dots, C_s$ are connected components of $T'$ as well. For all $i=1,\dots, m$, let $\bar C_i=\{\{v_2, w_i\}\}\cup C_i$. Then 
\begin{equation}\label{triangle}T'=\{\{v, v_2\}\}\cup\left(\bigcup_{i=1}^m\bar C_i\right)\cup\left(\bigcup_{i=m+1}^s C_i\right),\end{equation}
\noindent
and
\begin{equation}\label{triangle2}T''=\bigcup_{i=1}^sC_i,\end{equation}
\noindent
where both  unions are disjoint. The induction hypothesis applies to $T''$ and consequently to $C_i$ for all $i=1\,\dots, s$; hence, in view of Proposition \ref{components}, (\ref{triangle2}) implies:
\begin{equation}\label{14star} 
A''=\,{\rm ara}\,I(T'')={\rm pd}\,I(T'')=\sum_{i=1}^s{\rm pd}\,I(C_i)=\sum_{i=1}^s{\rm ara}\,I(C_i).\end{equation}
\noindent
For all $i=1,\dots, m$, let $A_i={\rm pd}\,I(\bar C_i)$. Since the number of vertices of $\bar C_i$ is less than the one of $T$, induction applies to $\bar C_i$, so that $A_i=\,{\rm ara}\,I(\bar C_i)$. We distinguish between two more subcases.\par\smallskip\noindent
\underline{Case 2.2.1}: We have that ara\,$I(C_i)\geq A_i$ for all $i=1,\dots, m$. Then, by (\ref{14star}), it holds:
\begin{eqnarray}\label{dot} 
A''&=&\sum_{i=1}^s{\rm ara}\,I(C_i)\nonumber\\
&\geq&\sum_{i=1}^m A_i + \sum_{i=m+1}^s{\rm ara}\, I(C_i)=\sum_{i=1}^m {\rm ara}\,I(\bar C_i)+\sum_{i=m+1}^s{\rm ara}\,I(C_i)\nonumber\\
&\geq&\,{\rm ara}\,\left(\sum_{i=1}^m I(\bar C_i)+\sum_{i=m+1}^sI(C_i)\right)=\,{\rm ara}\,I(T'\setminus\{\{v, v_2\}\})\geq A'-1. 
\end{eqnarray}
\noindent
  Recall from (\ref{q}) that $q'_{k_1}=v_2w_1+x_1y_1$. By the first equality of (\ref{T}) and (\ref{dot}) we deduce that
$$\Sigma: q_0=vv_2,\ q_1=vv_1+v_2w_1,\ q_2=x_1y_1,\ \Sigma'\setminus\{q'_0, q'_{k_1}\}$$
\noindent
is a tree-like system for $I(T)$ of length $A'+1\leq A''+2\leq A$ fulfilling (i). Hence the auxiliary claim is true in Case 2.2.1.\par\smallskip\noindent
 Before discussing Case 2.2.2, we need to show that there are subtrees $\Sigma_1,\dots,\Sigma_m$ of $\Sigma'$ such that $\Sigma_i$ is a tree-like system for $I(\bar C_i)$ for all $i=1,\dots, m$. This will be done in two steps. 
\par\bigskip\noindent
\underline{Claim 1}: Suppose that, for some element $q$ of $\Sigma'$, $q=a+b$, where $a$ and $b$ are edge monomials of  the edges $\alpha$ and $\beta$, respectively, and $\alpha\in \bar C_i$  for some $1\in\{1,\dots, m\}$. Then $\beta\in \bar C_i$.\newline
\underline{Proof of Claim 1}: For all $j=m+1,\dots, s$, $\bar C_i$ and $C_j$ are not connected to each other; therefore, as a consequence of Lemma \ref{path},  $\beta\not\in C_j$. Moreover, $\beta\ne\{v, v_2\}$, since, by the assumption of Case 2, $vv_2$ is an isolated summand in $\Sigma'$.  Consequently, in view of (\ref{triangle}) we have that $\beta\in \bar C_j$ for some $j\in\{1,\dots, m\}$. We show that $j=i$.  By Remark \ref{remark2}, there is an edge $\gamma$ lying between $\alpha$ and $\beta$. Then $\gamma\not\in C_j$ for all $j=m+1,\dots, s$. Moreover, since $\{v, v_2\}$ is a terminal edge of $T'$, we have that $\gamma\ne\{v, v_2\}$.  Hence $\gamma\in \bar C_h$ for some $h\in\{1,\dots,m\}$. Suppose that $h\ne i$. Then the common endpoint of $\alpha$ and $\gamma$ is the only common vertex of $\bar C_i$ and $\bar C_h$, namely $v_2$. But then $\alpha=\{v_2, w_i\}$, so that $q=q'_{k_i}$; therefore $\gamma=\gamma_i$, so that $v_2\in\gamma_i$, which contradicts the assumption of Case 2.2.. Hence $\gamma\in \bar C_i$. The arguments used for $\alpha$ can be applied to $\beta$,  which allows us to conclude that $j=i$, as required.
\par\bigskip\noindent
\underline{Claim 2}: For all $i=1,\dots, m$, let $\Sigma_i$ be the set of elements of $\Sigma'$ whose summands are edge monomials of edges of $\bar C_i$. Then $\Sigma_i$ is a subtree of $\Sigma'$ (and, consequently, a tree-like system for $I(\bar C_i)$).  \newline
\underline{Proof of Claim 2}: Let $\alpha\in \bar C_i$ and let $q$ be the element of $\Sigma'$ containing the edge monomial $a$ of $\alpha$. By Lemma \ref{strict} there is a unique strict subtree $\Sigma_{\alpha}$ of $\Sigma'$   ending at $q$. We show that $\Sigma_{\alpha}\subset \Sigma_i$. It will follow that 
$$\Sigma_i=\bigcup_{\alpha\in \bar C_i}\Sigma_{\alpha},$$
\noindent 
which, in view of Remark \ref{remark'}, will imply the claim. The above inclusion is obvious if $q=a$ is an isolated summand, because then $\Sigma_{\alpha}$ coincides with $q$.  So assume that $q=a+b$, where $b$ is the edge monomial of the edge $\beta$. From Claim 1 we know that $\beta\in \bar C_i$; from its proof we also know that the precedessor of $q$ in $\Sigma_{\alpha}$ contains an edge monomial $c$ whose corresponding edge $\gamma$ belongs to $\bar C_i$. This allows us to conclude by finite induction that all summands of $\Sigma_{\alpha}$ correspond to edges of $\bar C_i$, as required. 
\par\smallskip\noindent 
We are now ready to complete the discussion of Case 2.2. 
\par\bigskip\noindent
\underline{Case 2.2.2}: We have that ara\,$I(C_i)<A_i$ for some $i\in\{1,\dots, m\}$. Then by induction there is a tree-like system $\bar\Sigma_i$ for $I(C_i)$ having length less than $A_i$. On the other hand, $I(\bar C_i)=I(C_i)+(v_2w_i)$. Therefore, $\tilde\Sigma_i: v_2w_i, \bar\Sigma_i$ is a tree-like system for $I(\bar C_i)$ having length at most $A_i$. But, as a consequence of Claim 2, $A_i\leq\lambda(\Sigma_i)$.  Hence, replacing $\Sigma_i$ by $\tilde\Sigma_i$ in $\Sigma'$ as described in Remark \ref{remark01} produces a tree-like system for $I(T')$ whose length is not greater than $A'$ (thus it is, necessarily, equal to $A'$) and which contains $v_2w_i$ as an isolated summand. This takes us back to Case 1. Hence the auxiliary claim is true in Case 2.2.2.   
\par\smallskip\noindent
We have thus proven that the auxiliary claim is true in Case 2.2, hence it is true in Case 2.
\par\bigskip\noindent
\underline{Case 3}: The edge monomial $vv_2$ does not appear as an isolated summand in $\Sigma'$. Then there is an element of $\Sigma'$ of the form $q'=vv_2+b$, for some edge monomial $b$ of an edge $\beta=\{u,z\}$ of $T'$. According to Remark \ref{remark2}, there is an edge $\gamma$ of $T'$ lying between $\{v, v_2\}$ and $\beta$. Since $v$ is a terminal vertex of $T'$,  we  conclude that $v\not\in\gamma$, so that $v_2\in\gamma$, whence $\gamma=\{v_2, w_i\}$ for some index $i\in\{1,\dots, m\}$. Therefore $w_i\in\beta$, say $u=w_i$. Moreover, we may assume that $i=1$,  so that $q'=vv_2+w_1z$, and $\gamma=\{v_2, w_1\}$.  Hence the precedessor of $q'$ in $\Sigma'$  is either $v_2w_1$ (if this edge monomial appears as an isolated summand in $\Sigma'$), or the element $q'_{k_1}$. In the latter case, after applying tree-inversion  to the strict subtree of $\Sigma'$ ending at $q'$,  the starting point of this subtree becomes $v_2w_1$. This element can be pushed to the top of $\Sigma'$: hence we will henceforth work under the assumption that the starting point of $\Sigma'$ is $q'_0=v_2w_1$. Note that after the tree-inversion, the element of $\Sigma'$ containing the summand $vv_2$ is 
\begin{equation}\label{q'}q'=vv_2+w_hz'\qquad\mbox{ for some index $h\in\{1,\dots, m\}$ and some neighbour $z'$ of $w_h$;}\end{equation}
\noindent
 moreover, the precedessor of $q'$ is $q_{k_h}$ if $h\ne 1$, or $q'_0$ if $h=1$.\newline
Once again, we have to distinguish between two subcases.\par\smallskip\noindent 
\underline{Case 3.1}: One of the following conditions holds. Either
\begin{list}{}{}
\item{(a)} $A'\leq\,A''+1$, or
\item{(b)} the edge monomial $v_2w_j$ appears in $\Sigma'$ as an isolated summand for some index $j\in\{2,\dots, m\}$.
\end{list}
\noindent
First suppose that (a) holds. In this case replace $q'$ by the tree-like system formed by the two isolated summands $vv_2$ and $w_hz'$ and push $vv_2$ to the top. Then replace the former starting point $v_2w_1$ of $\Sigma'$ with $vv_1+v_2w_1$ and push it to the second position. This produces a tree-like system $\Sigma$  for $I(T)$ that has length $A'+1\leq A''+2\leq A$  and fulfills the auxiliary claim in case  (i). 
\newline
Now assume that (b) holds. We interchange $vv_2$ and $v_2w_j$ in $\Sigma'$. After this operation, $vv_2$ becomes an isolated summand  and $q'$ is turned into $q'=v_2w_j+w_hz'$.  The terminal edge $\{v, v_2\}$ of $T'$ cannot lie  between two edges, hence this operation does not affect the followers of $q'$. Moreover, the precedessor of $q'$ is still $q_{k_h}$ or $q'_0$. Hence, after this modification, $\Sigma'$ remains a tree-like system for $I(T')$; since it contains both $vv_2$ and $v_2w_1$ as isolated summands, we are taken back to Case 1. Hence the auxiliary claim is true in Case 3.1.
\par\bigskip\noindent
The next claim, which will be useful in the discussion of Case 3.2, is a consequence of Corollary \ref{corollary_connected}.
\par\smallskip\noindent
\underline{Claim 3}: For all $i=m+1,\dots, s$, let $\Sigma_i$ be the set of elements of $\Sigma'$ that contain the edge monomial of an edge of $C_i$. Then $\Sigma_i$ is a subtree of $\Sigma'$ and a tree-like system for $I(C_i)$. 
\par\bigskip\noindent
As a consequence, $\lambda(\Sigma_i)\geq\,{\rm ara}\,I(C_i)$ for all $i=m+1\,\dots, s$. Note that, in fact, 
\begin{equation}\label{asterisk} \lambda(\Sigma_i)=\,{\rm ara}\,I(C_i),\qquad(i=m+1,\dots, s),\end{equation}
\noindent
 because, otherwise, if we had $\lambda(\Sigma_i)>\,{\rm ara}\,I(C_i)$ for some index $i\in\{m+1,\dots, s\}$, by induction we could replace the subtree $\Sigma_i$ with an equivalent tree-like system $\Sigma'_i$ of length less than $\Sigma_i$; this would produce a tree-like system equivalent to $\Sigma'$ and of length smaller than $A'$, which is impossible.
\par\bigskip\noindent
\underline{Case 3.2}: We have that $A'\geq\,A''+2$ and none of the edge monomials $v_2w_i$ appears in $\Sigma'$ as an isolated summand for  $i\in\{2,\dots, m\}$. 
First assume that $m=1$, so that, according to the second equality of (\ref{T}), $I(T')=I(T'')+(vv_2, v_2w_1)$. It follows that $A'\leq A''+2$, so that, in view of the first part of the current assumption, $A'=A''+2$. Thus $\tilde\Sigma: vv_2, v_2w_1, \Sigma''$ is a tree-like system for $I(T')$ that has length $A'$ and to which Case 1 applies. So suppose that $m\geq 2$. 
Then, by the second part of the current assumption, for all $i=2,\dots, m$, $\Sigma'$ contains the element $q'_{k_i}$ defined in (\ref{q}). Moreover, $v_2$ has degree greater than 2; since $T$ is stretched, it follows that each of its neighbours $w_i$ has at most one neighbour $x_i$ other than $v_2$. In particular, in (\ref{q'}) we have that $z'=x_h$, so that $q'=vv_2+w_hx_h$. Finally, for all $i=2,\dots, m$, we have $q'_{k_i}=v_2w_i+d_i$, where $d_i$ is the edge monomial of the edge $\delta_i$ such that 
$$\delta_i=\begin{cases}{
 \{w_{j_i}, x_{j_i}\}\in C_{j_i}&  for some $j_i\in\{1,\dots, m\}$, $j_i\ne i$, if $i$ fulfills (\ref{q1});\cr\cr
\{x_i, y_i\}\in C_i& if $i$ fulfills (\ref{q2}).}
\end{cases}
$$
\noindent 
Note that no index $j_i$ can be equal to $h$, since otherwise both $q'_{k_i}$ and $q'$ would contain the edge monomial $w_hx_h$, against the definition of tree-like system. Moreover, case (\ref{q2}) cannot occur for $i=h$, since otherwise each of $q'_{k_h}=v_2w_h+x_hy_h$ and $q'$ would be the precedessor of the other in $\Sigma'$: in fact $v_2w_h$ divides $vv_2\cdot w_hx_h$ and $w_hx_h$ divides $v_2w_h\cdot x_hy_h$. But this is clearly a contradiction.  Hence we can define the following map:
$$\phi:\{2,\dots, m\}\longrightarrow \{1,\dots,m\}\setminus\{h\}$$
$$i\mapsto \begin{cases}{
j_i & if  $i$ fulfills (\ref{q1});\cr\cr
i & if  $i$ fulfills (\ref{q2}).}
\end{cases}
$$
\noindent
We show that the map $\phi$ is bijective. It suffices to prove injectivity.   Suppose for a contradiction that  we have $\phi(i)=\phi(i')$ for some $i, i'\in\{2,\dots, m\}$, $i\ne i'$. Since $i$ and $i'$ are distinct, (\ref{q2}) cannot hold for $i$ and $i'$ at the same time. First suppose that (\ref{q1}) occurs for $i$ and (\ref{q2}) occurs for $i'$. Then $i'=j_i$ and we have
\begin{eqnarray*} q'_{k_i}&=&v_2w_i+w_{j_i}x_{j_i}=v_2w_i+w_{i'}x_{i'}\\
q'_{k_{i'}}&=&v_2w_{i'}+x_{i'}y_{i'}.
\end{eqnarray*}
\noindent
It follows that $q'_{k_i}$ and $q'_{k_{i'}}$ are each one the precedessor of the other, which is impossible. So suppose that (\ref{q1}) occurs for $i$ and $i'$. Then $j_i=j_{i'}$, whence
 \begin{eqnarray*} q'_{k_i}&=&v_2w_i+w_{j_i}x_{j_i}\\
q'_{k_{i'}}&=&v_2w_{i'}+w_{j_{i'}}x_{j_{i'}}=v_2w_{i'}+w_{j_i}x_{j_i}.
\end{eqnarray*}
\noindent
This again is impossible, since $q'_{k_i}$ and $q'_{k_{i'}}$ both contain the edge monomial $w_{j_i}x_{j_i}$, against the definition of tree-like system.  This shows that $\phi$ is bijective. Note that its surjectivity implies  that for all $i=1,\dots, m$, $w_i$ has exactly one neighbour $x_i\ne v_2$;  we already knew it for $i=h$, in view of the form of $q'$. For all $j=1,\dots, m$, let $\bar S_j$ be the set of all elements of $\Sigma'$ which contain the edge monomial of  an edge of $C_j$, and let $S_j$ be the subset of those elements whose edge monomials all fulfill this condition. We investigate the relation between $\bar S_j$ and $S_j$. Let $q\in \bar S_j\setminus S_j$; then $q=a_j+b_j$ where $a_j$ and $b_j$ are the edge monomials of some edges $\alpha_j$ and $\beta_j$, such that, up to exchanging summands, $\alpha_j\in C_j$, $\beta_j\notin C_j$. Then $\beta_j\ne \{v_2, w_1\}$, because we are assuming that $v_2w_1$ is an isolated summand. Furthermore, for all indices $i\in\{1,\dots, m\}$, $i\ne j$, we have that $\beta_j\notin C_i$: otherwise the edge lying between $\alpha_j$ and $\beta_j$ would have an endpoint in $C_j$ and the other one in $C_i$, which is impossible. In view of (\ref{triangle}), there are only the following cases left: either $\beta_j=\{v, v_2\}$, or $\beta_j\in\bar C_i\setminus C_i$ for some index $i\in\{2,\dots,m\}$. In the former case, $q=q'$, i.e., $a_j=w_hx_h$, so that $j=h$. In the latter case, $b_j=v_2w_i$, so that $j=\phi(i)$, whence $j\ne h$, and  $q=q_{k_i}$, i.e., $a_j=d_i$. Thus
\begin{eqnarray}\label{minus}
\bar S_j&=&S_j\cup\{q_{k_{\phi^{-1}(j)}}\}\qquad(j\in\{1,\dots, m\}\setminus\{h\}),\nonumber\\
\bar S_h&=&S_h\cup\{q'\}.
\end{eqnarray}
\noindent
The set of summands of $\bar S_j$ is equal to the set of edge monomials of $I(C_j)$ together with $v_2w_{\phi^{-1}(j)}$ for all $j\in\{1,\dots, m\}\setminus\{h\}$, whereas the set of summands of $\bar S_h$ is equal to the set of edge monomials of $I(C_h)$ together with $vv_2$. In view of (\ref{triangle}) and Claim 3 it follows that
\begin{equation}\label{sigma}\Sigma'=\bigcup_{i=1}^m\bar S_i\,\cup\, \bigcup_{i=m+1}^s\Sigma_i\,\cup\,\{v_2w_1\},\end{equation}
\noindent
where the union is disjoint.  Set
\begin{equation}\label{plus}\Sigma'_i=S_i\cup\{a_i\}\qquad(i\in\{1,\dots, m\}).\end{equation}
\noindent
In other words, for all $i=1,\dots, m$, $\Sigma'_i$ is obtained from $\bar S_i$ by omitting the only edge monomial (namely, $b_i$) which does not belong to $I(C_i)$. Since, by Lemma \ref{path}, the edge lying between two edges of $C_i$ belongs to $C_i$ as well, it easily follows that $\Sigma'_i$ is a tree-like system for $I(C_i)$. For all $i=1,\dots, m$, set $\lambda_i=\lambda(\Sigma'_i)$. Then, comparing (\ref{minus}) and (\ref{plus}) we deduce that  $\lambda_i=\vert\bar S_i\vert$. Therefore, according to (\ref{14star}), the assumption of Case 3.2, (\ref{sigma}) and (\ref{asterisk})   we have that
\begin{eqnarray*}
\sum_{i=1}^m{\rm ara}\,I(C_i)+\sum_{i=m+1}^s{\rm ara}\,I(C_i)+2&=&A''+2\leq A'=\lambda(\Sigma')\\
&=&\sum_{i=1}^m\lambda_i+\sum_{i=m+1}^s\lambda(\Sigma_i)\,+1\\
&=&\sum_{i=1}^m\lambda_i+\sum_{i=m+1}^s{\rm ara}\,I(C_i)+1.\\
\end{eqnarray*}
\noindent
This implies that 
\begin{equation}\label{condition1}  {\rm ara}\,I(C_i)<\lambda_i\qquad\mbox{ for some index }i\in\{1,\dots, m\}.\end{equation}
\noindent 
By induction, we deduce that there is a tree-like system $\bar\Sigma_i$ of length less than $\lambda _i$ for $I(C_i)$. In $\Sigma'$ replace $\bar S_i$ with the tree-like system $\tilde\Sigma_i: v_2w_{\phi^{-1}(i)}, \bar\Sigma_i$, if $i\in\{1,\dots, m\}\setminus\{h\}$, or $\bar S_h$ with the tree-like system: $\tilde\Sigma_h: vv_2, \bar\Sigma_h$ if $i=h$, as described in Remark \ref{remark01}. In both cases $\Sigma'$ is replaced by an equivalent tree-like system of non greater length (hence, of the same length), where either $v_2w_j$ for some $j\ne 1$ or $vv_2$ is an isolated summand. We are thus taken back to Case 3.1 (b) or to Case 1 respectively. 
This shows that the auxiliary claim is true in Case 3.2., hence it is true in Case 3.  
\newline
We have thus shown the theorem for $n=2$. In order to complete the induction step, suppose now that $n\geq3$, $m=1$, $N\geq 4$, that the theorem is true for all smaller values of $N$, and that the auxiliary claim is true for all smaller values of $n$. 
 Consider the subgraphs $T'$ and $T''$ of $T$ introduced above. There is a tree-like system $\Sigma'': q''_0,\dots, q''_{A''-1}$ for  $I(T'')$.   The neighbours of $v$ in $T'$ are the $n-1$ vertices $v_1,\dots, v_{n-2}, v_n$, with $v_1, \dots, v_{n-2}$ of degree 1. Hence, by induction, there is a tree-like system $\Sigma': q'_0,\dots q'_{A'-1}$ for $I(T')$ fulfilling the auxiliary claim, i.e., such that
\begin{list}{}{}
\item{(i)$'$}  if $n=3$, then
\begin{eqnarray*}
q'_0&=&vv_3,\\
q'_1&=&vv_1+v_3w_1;
\end{eqnarray*}
\noindent
\item{(ii)$'$} if $n\geq4$, then 
\begin{eqnarray*}
q'_0&=&vv_n,\\
q'_1&=&vv_1+v_nw_1,\\
q'_2&=&vv_2,\\
&\vdots&\\
q'_{n-2}&=&vv_{n-2}.
\end{eqnarray*}
\end{list}
\noindent
First assume that $n=3$. Then  
\begin{equation}\label{rad'}\sqrt{(q'_0, q'_1)}=(vv_3, vv_1, v_3w_1).\end{equation}
\noindent
Set $q_2=vv_2$, and $q_i=q'_i$ for $i=0,1$. From (\ref{rad'}) it follows that
$$\sqrt{(q_0, q_1, q_2)}=\sqrt{(q'_0, q'_1, vv_2 )}=(vv_3, vv_1, v_3w_1, vv_2).$$
\noindent
Hence
\begin{eqnarray*}
I(T)&=&(vv_3, vv_1, vv_2,  v_3w_1)+I(T'')\\
&=&\sqrt{(q_0, q_1, q_2)}+\sqrt{(q''_0,\dots, q''_{A''-1})},
\end{eqnarray*}
\noindent
which shows that, if $n=3$, $I(T)$ is generated, up to radical, by the following tree-like system of length $A''+3\leq A$:
\begin{equation}\label{elements}q_0, q_1, q_2, q''_0,\dots, q''_{A''-1}.\end{equation}
\noindent
Moreover, (\ref{elements}) together with (i)$'$, tells us that (i) is  true for $T$.\newline
Now assume that $n\geq4$.   We have: 
\begin{equation}\label{rad}\sqrt{(q'_0,\dots, q'_{n-2})}=(vv_n, vv_1, \dots, vv_{n-2}, v_nw_1).\end{equation}
\noindent
Set $q_{n-1}=vv_{n-1}$, and $q_i=q'_i$ for all indices $i=0,\dots, n-2$. From (\ref{rad}) it follows that
$$\sqrt{(q_0, \dots, q_{n-2}, q_{n-1})}=\sqrt{(q'_0, \dots, q'_{n-2}, vv_{n-1} )}=(vv_n, vv_1, \dots, vv_{n-2}, v_nw_1, vv_{n-1}).$$
\noindent
Hence 
\begin{eqnarray*}
I(T)&=&(vv_n, vv_1, \dots, vv_{n-2}, vv_{n-1},  v_nw_1)+I(T'')\\
&=&\sqrt{(q_0,\dots, q_{n-2}, q_{n-1})}+\sqrt{(q''_0,\dots, q''_{A''-1})},
\end{eqnarray*}
\noindent
which shows that, if $n\geq4$, $I(T)$ is generated, up to radical, by the following tree-like system  of length $A''+n\leq A$:
\begin{equation}\label{elements'}q_0,\dots, q_{n-2}, q_{n-1}, q''_0,\dots, q''_{A''-1}.\end{equation}
\noindent
Moreover, (\ref{elements'}) together with (ii)$'$ tells us that (ii) is true for $T$. 
This completes the proof of the claim. 
\par\bigskip\noindent
The proof of the theorem we have given is in fact a constructive one. Following the thread of arguments developed there,  one can recursively produce, for any stretched forest $T$,  a tree-like system of  length ara\,$I(T)$ for $I(T)$, by induction on the number of vertices. In all the possible cases we have described how to (easily) obtain the required tree-like system for $I(T)$ from a tree-like system for $I(T')$ or $I(T'')$. We explicitly exploited the assumption that the forest $T$ is stretched; this, evidently, played a crucial in the treatment of Case 3.2. We do not know how to remove it, in order to extend the procedure to all forests.  In fact, the next result reveals the peculiarity of stretched forests: the associated tree-like systems share a property which does not hold, in general, for non stretched ones. This result makes Lemma \ref{strict} more precise. 
\begin{proposition}\label{corollary3} If $T$ is a stretched forest, then all tree-like systems whose support is the set of edge monomials of $I(T)$ are the disjoint union of strict subtrees. 
\end{proposition}
\demo Let $\Sigma$ be a tree-like system for $I(T)$, and let $q$ be an element of $\Sigma$. From Lemma \ref{strict} we know that there is a subtree of $\Sigma$ containing $q$. We show that this subtree is unique. Since, by Lemma \ref{strict}, there is a unique strict subtree of $\Sigma$ ending at $q$, it suffices to show that $q$ cannot have more than one follower. First assume that $q=a$ is an isolated summand, where $a$ is the edge monomial of  the edge $\alpha=\{u,v\}$ of $T$. If $q$ had two different followers, then these would be  $q_1=uu_1+vv_1$ and $q_2=uu_2+vv_2$, where $u_1, u_2$ and $v_1, v_2$ are distinct neighbours of $u$ and $v$ respectively. But then $u$ and $v$  would both have degree at least 3, against the definition of stretched forest.  Now assume that $q=a+b$, where $a$ is as above and $b$ is the edge monomial of the edge $\beta=\{w,z\}$ of $T$. Suppose for a contradiction  that $q$ has two different followers $q_1=c_1+d_1$ and $q_2=c_2+d_2$, where $c_i$ and $d_i$ are the edge monomials of the edges $\gamma_i$ and $\delta_i$ respectively, for $i=1,2$.   If the edge lying between $\gamma_i$ and $\delta_i$ is $\alpha$ for $i=1,2$, then $q_1$ and $q_2$ are of the same form as in the previous case, which, as we have seen, leads to a contradiction. We come to the same conclusion if the edge lying between $\gamma_i$ and $\delta_i$ is $\beta$ for $i=1,2$. So we have to assume that, up to exchanging indices, $\alpha$ lies between $\gamma_1$ and $\delta_1$ and $\beta$ lies between $\gamma_2$ and $\delta_2$. Then $q_1=uu_1+vv_1$ and $q_2=ww_1+zz_1$, where $u_1$, $v_1$, $w_1$ and $z_1$ are neighbours of $u, v, w$ and $z$, respectively. Up to renaming vertices, we may assume that the edge lying between $\alpha$ and $\beta$ is $\gamma=\{v, w\}$.  We conclude that  $v$ and $w$ both have degree at least 3: the vertices $u, 
v_1, w$ are neighbours of $v$, the vertices $z, w_1, v$ are neighbours of $w$.   This, once again, contradicts the definition of stretched forest. This completes the proof. 
\par\bigskip\noindent
\begin{remark}{\rm We conjecture that the theorem is true for all forests. Some of the examples contained in the next section will provide supporting evidence for this.\newline
In fact, every forest can be obtained from a stretched one by replacing some of the subgraphs $\{\{u, u'\},\{u', u''\}\}$ (where $u'$ has degree 2) with the edge $\{u, u''\}$. We, however, cannot predict, in general, in which way this operation affects the projective dimension or the arithmetical rank. }\end{remark}
\begin{example}\label{Iexample}{\rm  Let $T$ be the stretched tree whose edges are
$$\{v, v_1\},\ \{v, v_2\},\ \{v, v_3\},\ \{v_3, w_1\},\  \{w_1, a\},\ \{w_1, b\},\ \{w_1, c\}.$$
\noindent
Here $n=3$, $m=1$. We have that pd\,$I(T)=\,$ara\,$I(T)=6$ and a tree-like system for $I(T)$  is:
$$\Sigma: q_0=vv_3,\ q_1=vv_1+v_3w_1,\ q_2=vv_2,\ q_3=w_1a,\ q_4=w_1b,\ q_5=w_1c.$$
\noindent
It fulfills the auxiliary claim, case (ii), of the proof of Theorem \ref{theorem}. Moreover, it is the disjoint union of five strict subtrees, the first one is formed by $q_0, q_1$, the other four are formed by the isolated summands $q_2, q_3, q_4, q_5$. 
}
\end{example}
\begin{example}\label{IIexample}{\rm  Consider the stretched tree whose edges are
$$\{v, v_1\},\ \{v, v_2\},\ \{v_2, w_1\},\ \{v_2, w_2\},\  \{w_1, a\},\ \{a, b\},\ \{a, c\},\ \{c, d\}.$$
\noindent
Here $n=2$, $m=2$. We have that pd\,$I(T)=\,$ara\,$I(T)=5$ and a tree-like system for $I(T)$ fulfilling the auxiliary claim is:
$$\Sigma: q_0=vv_2,\ q_1=vv_1+v_2w_1,\ q_2=v_2w_2+w_1a,\ q_3=ac,\ q_4=ab+cd.$$
\noindent
It fulfills the auxiliary claim, case (i), of the proof of Theorem \ref{theorem}. It is the disjoint union of two strict subtrees, formed by $q_0, q_1, q_2$ and by $q_3, q_4$ respectively. 

}
\end{example}
\section{The arithmetical rank of some special trees}
In this section we present some special classes of trees for which the arithmetical rank equals the projective dimension. In all cases, we will show that there is a tree-like system which generates the edge ideal up to radical, even if the tree is not stretched. \subsection{Double-star graphs}
Let $r, s$ be non negative integers, and consider the tree on the vertex set $\{a, b, x_1,\dots, x_r, y_1,\dots, y_s\}$ given as follows:
$$T_{r,s}=\{\{a,b\}\}\cup\{\{a, x_i\}:i=1,\dots, r\}\cup\{\{b,y_i\}:i=1,\dots, s\}.$$
\noindent
 We will call $T_{r,s}$ a {\it double-star} graph. 
 We consider the edge ideal $I(T_{r,s})$ of $T_{r,s}$ in the polynomial ring $$R=K[a,b, x_1,\dots, x_r, y_1,\dots, y_s].$$ The edge monomials of of $I(T_{r,s})$ are the following $r+s+1$ elements:
\begin{equation}\label{generators}
ab,\ ax_1,\ \dots, \ ax_r,\ by_1,\ \dots, \ by_s.
\end{equation}
\noindent
According to \cite{J}, Example 2.1.7, 
$${\rm pd}\,I(T_{r,s})=\max\{r,s\}+1,$$
\noindent
so that
\begin{equation}\label{ara2}{\rm ara}\,I(T_{r,s})\geq\max\{r,s\}+1.\end{equation}
\noindent
In $T_{r,s}$ all vertices $x_1,\dots, x_r$ and $y_1,\dots, y_s$ have degree 1. Therefore, with respect to the notation introduced in Section 1, and in view of Remark \ref{remark}, 
$$\rho(I(T_{r,s}))=\min\{\deg(a), \deg(b)\}=\min\{r+1,s+1\}=\min\{r,s\}+1.$$ 
\noindent
Consequently,
\begin{eqnarray}\label{sharp}\mu(I(T_{r,s}))-\rho(I(T_{r,s}))+1 &=&r+s+1-\min\{r,s\}-1+1\nonumber\\
&=&r+s+1-\min\{r,s\}\nonumber\\
&=&\max\{r,s\}+1.
\end{eqnarray}
\noindent
From (\ref{ara2}) and (\ref{sharp}) we see that, for the edge ideal of any double-star graph, equality holds in (\ref{inequalities2}), i.e., the upper bound given in \cite{B}, Proposition 1 is sharp. Using Lemma \ref{lemma} one can easily prove that 
\begin{list}{}{}
\item{-} if $r\leq s$, then $I(T_{r,s})$ is generated up to radical by
$$ab,\ ax_1+by_1,\ \dots, \ ax_r+by_r,\ by_{r+1},\ \dots, \ by_s;$$
\item{-}
if $r> s$, then $I(T_{r,s})$ is generated up to radical by
$$ab,\ ax_1+by_1,\ \dots, \ ax_s+by_s,\ ax_{s+1},\ \dots, \ ax_r.$$
\end{list}
\noindent
In both cases,  we have a tree-like system of length ara\,$I(T_{r,s})$ for $T_{r,s}$. However, $T_{r,s}$ is stretched if and only if $r=1$ or $s=1$. 
 If $s=0$, then $T_{r,s}$ is the {\it star-graph} $S_{r+1}$. 
We have
$${\rm ara}\,I(S_{r+1})={\rm pd}\,I(S_{r+1})=\mu(I(S_{r+1}))-\rho(I(S_{r+1}))+1=r+1.$$
\subsection{Line graphs}
The star graph $S_1$ is a graph with one edge only. It is the simplest example of {\it line graph}: for every integer $r\geq2$, consider the following stretched tree on the  vertex set $\{x_1, \dots, x_r\}$:
 $$L_r=\{\{x_i, x_{i+1}\}: i=1,\dots, r-1\}.$$
\noindent
Its edge ideal in $R=K[x_1, \dots, x_r]$ is 
$$I(L_r)=(x_1x_2,\ x_2x_3, \ \dots, \ x_{r-1}x_r).$$
\noindent
The projective dimension of $L_r$ has been completely determined in \cite{J}, Corollary 7.7.35. We have to distinguish between three cases, depending on the  residue of $r$ modulo 3. In each case, we make use of Lemma \ref{lemma} for determining elements generating ideal $I(L_r)$ up to radical. 
\begin{list}{}{}
\item{-} If $r=3s$ for some integer $s$, then ${\rm pd}\, I(L_r)={\rm ara}\, I(L_r)=2s$ and $I(L_r)$ is generated up to radical by $x_1x_2, x_2x_3$ if $s=1$ and, if $s\geq 2$, by:
\begin{eqnarray*}
&&x_2x_3,\ x_1x_2+x_3x_4,\\
&&\vdots\\
&&x_{3k-1}x_{3k},\ x_{3k-2}x_{3k-1}+x_{3k}x_{3k+1},\\
&&\vdots\\
&&x_{3s-4}x_{3s-3},\ x_{3s-5}x_{3s-4}+x_{3s-3}x_{3s-2},\\
&&x_{3s-2}x_{3s-1},\\
 &&x_{3s-1}x_{3s}.
\end{eqnarray*}
\item{-} If $r=3s+1$ for some integer $s$, then ${\rm pd}\, I(L_r)={\rm ara}\, I(L_n)=2s$ and $I(L_r)$ is generated up to radical by:
\begin{eqnarray*}
&&x_2x_3,\ x_1x_2+x_3x_4,\\
&&\vdots\\
&&x_{3k-1}x_{3k},\ x_{3k-2}x_{3k-1}+x_{3k}x_{3k+1},\\
&&\vdots\\
&&x_{3s-1}x_{3s},\ x_{3s-2}x_{3s-1}+x_{3s}x_{3s+1}.
\end{eqnarray*}
\item{-} If $4=3s+2$ for some integer $s$, then ${\rm pd}\, I(L_r)={\rm ara}\, I(L_r)=2s+1$ and $I(L_r)$ is generated up to radical by:
\begin{eqnarray*}
&&x_2x_3,\ x_1x_2+x_3x_4,\\
&&\vdots\\
&&x_{3k-1}x_{3k},\ x_{3k-2}x_{3k-1}+x_{3k}x_{3k+1},\\
&&\vdots\\
&&x_{3s-1}x_{3s},\ x_{3s-2}x_{3s-1}+x_{3s}x_{3s+1},\\
&&x_{3s+1}x_{3s+2}.
\end{eqnarray*}
\end{list}
\noindent
Ideal $I(L_3)=(x_1x_2, x_2x_3)$ is a complete intersection. For $r\geq 4$ the auxiliary claim in the proof of Theorem \ref{theorem} is fulfilled for $v=x_2$, $n=2$, $m=1$, $v_1=x_1$, $v_2=x_3$ and $w_1=x_4$.\newline 
Also note that, for all $r\geq 2$, $\mu(I(L_r))=r-1$, whereas $\rho(I(L_2))=\rho(I(L_3))=1$, and $\rho(I(L_r))=2$ for all $r\geq 4$. An elementary computation shows that equality holds in  (\ref{inequalities2}) for $I(L_r)$, i.e., ${\rm ara}\,I(L_r)=\mu(I(L_r))-\rho(I(L_r))+1$, if and only if $2\leq r\leq 6$. In all the other cases the inequality is strict.
\section{An addition on double-star graphs}
In this section we explicitly determine the minimal free resolution of the edge ideal of the double-star graph $T_{r,s}$ introduced above. Our approach is independent from the one adopted by Jacques and Katzman in \cite{J} and \cite{JK}. We follow the method developed by Lyubeznik in \cite{L1}. Let us recall how he explicitly constructs, for any monomial ideal, 
  a free resolution which is obtained from  the well-known 
  Taylor resolution by omitting redundant terms.\par\noindent 
Let $f_1, \dots,f_m$ be an {\sl ordered} sequence of $s$  monomials of the polynomial ring $R$ over a field, let $I$ be the ideal generated by these monomials.  
\begin{definition}\label{definition5} {\rm For all sequences 
$(f_{i_1},\dots,f_{i_t})$, where $1\leq i_1<\cdots<i_t\leq m$, the symbol 
$u(f_{i_1},\dots,f_{i_t})$ will be called {\it L-admissible of dimension $t$} if $f_q$ 
does not divide lcm$\,(f_{i_h},f_{i_{h+1}},\dots,f_{i_t})$ for all
$h<t$ and $q<i_h$.}
\end{definition}
 Set $L^0=R$ and for all $t=1,\dots,m$, let $L^{t}$ be the free $R$-module generated by all
$L$-admissible symbols of dimension $t$. Define the map $d_t:L^{t}\to L^{t-1}$
by setting
\begin{equation}\label{d}d_t(u(f_{i_1},\dots,f_{i_t}))=\displaystyle\sum_{j=1}^t(-1)^{j+1}\frac{{\rm lcm}
(f_{i_1},\dots,f_{i_t})}{{\rm lcm}(f_{i_1},\dots,\hat{f_{i_j}},\dots,
f_{i_t})}u(f_{i_1},\dots,\hat{f_{i_j}},\dots,f_{i_t}).\end{equation}
Then one has the following
\begin{theorem}\label{Lresolution}{\rm (\cite{L1}, p.~193)} The complex
$$\qquad\qquad\qquad\qquad\qquad\qquad 0\to L^{m}\buildrel d_m\over\to L^{m-1}\buildrel d_{m-1}\over\to\cdots
\buildrel d_1\over\to L^0\to 0\qquad\qquad\qquad\qquad(\ast)$$
is a free resolution of $R/I$.
\end{theorem}
The resolution $(\ast)$ is called a {\it Lyubeznik resolution} of $I$. 
Note that the Lyubeznik resolution of $I$ in general strictly depends on the
order of the sequence $f_1,\dots,f_m$: different permutations of the
$f_i$ can give rise to non-isomorphic resolutions.  Note that resolution ($\ast$) is minimal if and only if for all admissible symbols $u(f_{i_1},\dots,f_{i_t})$,
\begin{equation}\label{lcm}\lcm(f_{i_1},\dots,
f_{i_t})\neq\lcm(f_{i_1},\dots,\hat{f_{i_j}},\dots,
f_{i_t})\mbox{ for all }j=1,\dots, t,\end{equation}
\noindent which is true if and only if, among the monomials $f_{i_1},\dots,f_{i_t}$, none divides the least common multiple of the remaining.\par\smallskip\noindent
Next we show that the above construction yields a minimal free resolution of $I(T_{r,s})$. On the minimal monomial generating set of $I(T_{r,s})$ given in (\ref{generators}) we fix the following order:
\begin{equation}\label{order}f=ab,\ g_1=ax_1,\ \dots,\ g_r=ax_r,\ h_1=by_1,\ \dots,\ h_s=by_s.\end{equation}
\noindent
Note that the only divisibility relations between the generating monomials listed in (\ref{order}) are those deriving from:
$$f\vert g_i\cdot h_j.\qquad(1\leq i\leq r, 1\leq j\leq s).$$
Consequently, the admissible symbols in the sense of Definition \ref{definition5} are those which do not contain both an element $g_i$ and element $h_j$, i.e., those having one of the following forms:
\begin{equation}\label{admissible}u(f),\ u(f,g_{i_1}, \dots, g_{i_p}),\ u(f,h_{i_1}, \dots, h_{i_q}),\ u(g_{i_1}, \dots, g_{i_p}),\ u(h_{i_1}, \dots, h_{i_q}),\end{equation}
\noindent
with $1\leq p\leq r$, $1\leq q\leq s$. These all fulfill condition (\ref{lcm}), because each generator contains an indeterminate which does divide any of the remaining generators appearing in the same admissible symbol.  We have thus proven:
\begin{proposition}\label{resolution} The Lyubeznik resolution associated with (\ref{order}) is a minimal free resolution of $I(T_{r,s})$. 
\end{proposition}
\noindent It follows that, for all $t=1,\dots, r+s+1$, the $t$-th Betti number $\beta_t$ of $I(T_{r,s})$ is the rank of $L^t$, i.e., the number of admissible symbols of dimension $t$ listed in (\ref{admissible}). The projective dimension of $I(T_{r,s})$ is the maximum dimension of these admissible symbols. 
\begin{proposition}\label{betti} We have that pd\,$I(T_{r,s})=\max\{r,s\}+1$. Moreover, for all $t=1,\dots, {\rm pd}\,I(T_{r,s})$, 
$$\beta_t=\begin{cases}{
1+r+s & if $t=1$
\cr\cr
{{r+1}\choose{t}}+{s+1\choose t} &  if $2\leq t\leq {\rm pd}\,I(T_{r,s})$,}
\end{cases}$$
where we have set equal to zero all binomial coefficients ${a\choose b}$ with $a<b$.
\end{proposition} 
\demo In (\ref{admissible}), an admissible symbol of maximum dimension is $u(f, g_1,\dots, g_r)$ if $r\geq s$, or $u(f, h_1,\dots, h_s)$ if $s\geq r$. This shows the first part of the claim.\newline 
It is well known that $\beta_1$ is the number of minimal generators of $I(T_{r,s})$, which are listed in (\ref{generators}), i.e.,  $\beta_1=r+s+1$. Now let $2\leq t\leq {\rm pd}\,I(T_{r,s})$. In view of (\ref{admissible}), the admissible symbols of dimension $t$ are those having one of the following forms:
\begin{list}{}{}
\item{(i)} $u(f,g_{i_1}, \dots, g_{i_{t-1}})$ (there are ${r\choose{t-1}}$ of this kind);
\item{(ii)} $u(g_{i_1}, \dots, g_{i_{t}})$ (there are ${r\choose{t}}$ of this kind);
\item{(iii)} $u(f,h_{i_1}, \dots, h_{i_{t-1}})$ (there are ${s\choose{t-1}}$ of this kind);
\item{(iv)} $u(h_{i_1}, \dots, h_{i_{t}})$ (there are ${s\choose{t}}$ of this kind).
\end{list}
\noindent    
Hence
$$\beta_t={{r}\choose{t-1}}+{r\choose t}+{s\choose{t-1}}+{s\choose t}={{r+1}\choose{t}}+{{s+1}\choose t},$$
\noindent
as was to be shown. This completes the proof.
\par\medskip\noindent
The 0-th syzygies of $I(T_{r,s})$ are its minimal monomial generators listed in (\ref{generators}), which are all of degree 2.  From (\ref{d}) we can derive the degree of the higher syzygies of $I(T_{r,s})$. For $t=2,\dots, {\rm pd}\,I(T_{r,s})$, the entries of the $(t-1)$-th syzygy matrix are of one of the following forms:
\begin{eqnarray*}\frac{\lcm(f,g_{i_1}, \dots, g_{i_{t-1}})}{\lcm(g_{i_1}, \dots, g_{i_{t-1}})}&=&b, \\
(-1)^j\frac{\lcm(f,g_{i_1}, \dots, g_{i_{t-1}})}{\lcm(f, g_{i_1}, \dots, \hat{g_{i_j}},\dots, g_{i_{t-1}})}&=&(-1)^jx_{i_j},\\
(-1)^{j+1}\frac{\lcm(g_{i_1}, \dots, g_{i_{t}})}{\lcm(g_{i_1}, \dots, \hat{g_{i_j}},\dots, g_{i_{t}})}&=&(-1)^{j+1}x_{i_j}, \\
\frac{\lcm(f,h_{i_1}, \dots, h_{i_{t-1}})}{\lcm(h_{i_1}, \dots, h_{i_{t-1}})}&=&a, \\
(-1)^j\frac{\lcm(f,h_{i_1}, \dots, h_{i_{t-1}})}{\lcm(f, h_{i_1}, \dots, \hat{h_{i_j}},\dots, h_{i_{t-1}})}&=&(-1)^jy_{i_j},\\
(-1)^{j+1}\frac{\lcm(h_{i_1}, \dots, h_{i_{t}})}{\lcm(h_{i_1}, \dots, \hat{h_{i_j}},\dots, h_{i_{t}})}&=&(-1)^{j+1}y_{i_j}, \\
\end{eqnarray*} 
\noindent
which are all of degree 1. 
This shows the next result.
\begin{corollary}\label{linear} The edge ideal of a double-star graph has a 2-linear resolution. 
\end{corollary}
\begin{remark}{\rm According to the characterization given by Fr\"oberg \cite{F}, the above statement can be rephrased in purely combinatorial terms  by saying that $I(T_{r,s})$ is the Stanley-Reisner ideal of a simplicial complex which is the {\it clique complex} of a {\it chordal} graph (see \cite{EGHP2}, pp.~9--10 for the definitions of these terms and see \cite{EGHP2}, Theorem 2.1 for a generalization of this result). An equivalent algebraic geometric formulation is the following: the projective subvariety of ${\bf P}^{r+s+1}$ associated with $I(T_{r,s})$ is a {\it small scheme} in the sense of Eisenbud, Green, Hulek and Popescu \cite{EGHP}.}
\end{remark}
\begin{example}{\rm  We explicitly compute the minimal free-resolution  of the double-star graph $T_{2,3}$, whose edge ideal is
$$I(T_{r,s})=(ab,\ ax_1,\ ax_2,\ by_1,\ by_2,\ by_3).$$
\noindent
According to Proposition \ref{betti}, pd\,$I(T_{r,s})=4$, and the Betti numbers of $I(T_{r,s})$ are: 
$$\beta_1=6,\quad \beta_2=9,\quad \beta_3=5,\quad \beta_4=1.$$
\noindent
We give the $t$-th syzygy matrix for $t=1,2, 3$, with respect to the basis
$$u(ab), u(ax_1), u(ax_2), u(by_1), u(by_2), u(by_3)$$
\noindent of $L^1$, the basis
$$u(ab, ax_1), u(ab, ax_2), u(ax_1, ax_2), u(ab, by_1),  
u(ab, by_2),$$
$$u(ab, by_3), u(by_1, by_2), u(by_1, by_3), u(by_2, by_3)$$
\noindent
of $L^2$, the basis
$$u(ab, ax_1, ax_2), u(ab, by_1, by_2), u(ab, by_1, by_3), u(ab, by_2, by_3),$$ 
$$u(by_1, by_2, by_3)$$
\noindent of $L^3$, and the basis
$$u(ab, by_1, by_2, by_3)$$
\noindent
of $L^4$:

$$S_1=\left(\begin{array}{cccccc}
-x_1 & b & 0 & 0 & 0 & 0\\
-x_2 & 0 & b & 0 & 0 & 0\\
0 & -x_2 & x_1 & 0 & 0 & 0\\
-y_1 & 0 & 0 & a & 0 & 0\\
-y_2 & 0 & 0 & 0 & a & 0\\
-y_3 & 0 & 0 & 0 & 0 & a\\
0 & 0 & 0 & -y_2 & y_1 & 0\\
0 & 0 & 0 & -y_3 & 0 & y_1\\
0 & 0 & 0 & 0 & -y_3 & y_2\\
\end{array}\right)
$$
\vskip.3truecm\noindent
$$S_2=\left(\begin{array}{ccccccccc}
x_2& -x_1 & b  & 0    & 0    & 0     &  0  & 0    & 0\\
0  & 0    & 0    & y_2  & -y_1 & 0     & a & 0    &0\\
0  & 0    & 0    & y_3  & 0    & -y_1  & 0   & a  &0\\
0  & 0    & 0    & 0    & y_3  & -y_2  & 0   & 0    &a\\
0  & 0    & 0    & 0    & 0    & 0     & y_3 & -y_2 &y_1\\
\end{array}\right)
$$
\vskip.3truecm\noindent
$$S_3=\left(\begin{array}{ccccc}
0&-y_3 & y_2 & -y_1 & a
\end{array}\right)
$$
}\end{example}

\begin{center}
{\sc Acknowledgements}
\end{center}
\noindent
The author is indebted to Mordechai Katzman and Adam van Tuyl for valuable suggestions concerning this paper.
\end{document}